\documentclass[a4paper,10pt]{article}

\title{Totally odd immersions of complete graphs in graph products}

\usepackage{amsmath,amssymb,amsthm}
\usepackage{lmodern}
\usepackage{float}
\usepackage{todonotes}{}
\usepackage[T1]{fontenc}
\usepackage[utf8]{inputenc}
\usepackage{xcolor}
\usepackage{fullpage}
\usepackage{array}
\usepackage{soul}
\usepackage[mode=buildnew]{standalone}
\usepackage{color}
\usepackage{fancyvrb}
\usepackage{tikz}
\usepackage{wasysym}
\usepackage{authblk}
\usepackage{graphicx}
\usepackage{comment}
\usepackage{enumerate}

\newtheorem{theorem}{Theorem}[section]

\newtheorem{corollary}[theorem]{Corollary}

\newtheorem{conjecture}[theorem]{Conjecture}
\newtheorem{claim}[theorem]{Claim}
\newtheorem{question}[theorem]{Question}

\usepackage{ulem}

\newcommand{\toi}{\mathrm{toi}}
\newcommand{\im}{\mathrm{im}}

\newcommand{\p}{\mathcal{P}}

\author[1]{Henry Echeverr\'ia}
\author[1,2]{Andrea Jim\'enez}
\author[3]{Suchismita Mishra}
\author[1]{Daniel~A.~Quiroz}
\author[1]{Mauricio Y\'epez}

\affil[1]{\small Instituto de Ingenier\'ia Matem\'atica-CIMFAV, Universidad de Valpara\'iso, Chile.} %Email: \texttt{daniel.quiroz@uv.cl}}

\affil[2]{Millennium Nucleus for Social Data Science (SODAS), Santiago, Chile}

\affil[3]{Institute of Mathematical Science, Chennai, India}
\date{}

\begin{document}

\maketitle

\begin{abstract}

   For a graph $G$, let $\im(G)$ denote the maximum integer $t$ such that $G$ contains $K_t$ as an immersion. A recent paper of Collins, Heenehan, and McDonald (2023) studied the behaviour of this parameter under graph products, asking how large can $\im(G\ast H)$ be in terms of $\im(G)$ and $\im(H)$, when~$\ast$ is one of the four standard graph products. We consider a similar question for the parameter $\toi(G)$ which denotes the maximum integer $t$ such that $G$ contains $K_t$ as a totally odd immersion. As an application, we obtain that no minimum counterexample to the immersion-analogue of the Odd Hadwiger Conjecture can be obtained from the Cartesian, direct (tensor), lexicographic or strong product of graphs. 

\end{abstract}

\noindent \textbf{Keywords:} immersion, totally odd immersion, strong product, Cartesian product, direct product

\section{Introduction}

To \underline{split off} a pair of adjacent edges $uv,vw$ (with $u,v,w$ distinct) is to delete both these edges and to add the edge $uw$, if it is not already present. A graph $G$ is said to contain another graph~$H$ as an \underline{immersion} if $H$ can be obtained from $G$ by iteratively splitting off edges and deleting vertices and edges. In 1965, Nash-Williams conjectured~\cite{NashWilliams} that graphs are well-quasi-ordered by the immersion relation, and his conjecture was confirmed~\cite{Robertson_Seymour} by Robertson and Seymour in 2010. By $\im(G)$ we denote the maximum integer $t$ such that $G$ contains an immersion of $K_t$. This is a well-studied parameter. DeVos, McDonald, Mohar and Scheide~\cite{DeVosMMS}, and, independently, Wollan~\cite{Wollan} gave a structure theorem for graphs~$G$ with $\im(G)<t$. Liu~\cite{Liu1} provided a different structure theorem for excluded immersions. Moreover, the following conjecture, which relates $\im(G)$ to the chromatic number, $\chi(G)$, of $G$ has received considerable attention in recent years \cite{BoJiLiPaQuiSa24,BustamanteQSZ2021, DeVosDFMMS,DeVosKMO2010,dvovrak2018complete,GauthierLW2019,Liu}.

\begin{conjecture}[Abu-Khzam and Langston \cite{AbuLangston}]\label{conj:imm}
    For every graph $G$ we have $\chi(G)\le \im(G)$.
\end{conjecture}

In a recent paper, Collins, Heenehan and McDonald~\cite{collins2023clique} studied the behaviour of complete graph immersions when one of the four standard graph products $-$ Cartesian, direct (or tensor), lexicographic, and strong product $-$ are applied. They raised the following question.

\begin{question}[Collins, Heenehan and McDonald~\cite{collins2023clique}]\label{qus:im}
    Let $G$ and $H$ be graphs with $\im(G) = t$ and $\im(H) = r$, and let $\ast$ be any of the four standard graph products. Is $\im(G \ast H) \geq \im(K_t \ast K_r)$?
\end{question}

In the same paper, they provide a positive answer to Question~\ref{qus:im} for the lexicographic and Cartesian products. In particular, these results imply that no minimal counterexample to Conjecture~\ref{conj:imm} can be obtained through lexicographic or Cartesian product. In addition, they conjectured that a positive answer also holds for the direct product and strong products. Some evidence towards these conjectures is given in their paper and in a paper of Wu and Deng \cite{wu2025note}. 

Our contributions are inspired by this line of research and are concerned with a restriction of the immersion relation which we now introduce.

\subsection{Totally odd immersions}

In a way that is equivalent to the definition through split off, we can say that graph~$G$ contains another graph $H$ as an immersion if there exists an injective function $\varphi\colon V(H)\rightarrow V(G)$ such that:
\begin{enumerate}[(I)]
\item For every $uv\in E(H)$, there is a path in $G$, denoted $P_{uv}$, with endpoints $\varphi(u)$ and~$\varphi(v)$.
\item The paths in $\{P_{uv} \mid uv\in E(H) \}$ are pairwise edge-disjoint.
\end{enumerate}
The vertices in $\phi(V(H))$ are called the \underline{terminals} of the immersion. If the terminals are not allowed to appear as interior vertices on paths in $P_{uv}$, then $G$ is said to contain $H$ as a \underline{strong immersion}, and if, moreover, 
the paths are internally vertex-disjoint, then indeed $G$ contains a \underline{subdivision} of $H$. An immersion (or strong immersion) 
is said to be \underline{totally odd} if all the paths are odd.

Churchley~\cite{churchley2017odd} gave a structure theorem for graphs excluding $K_t$ as a totally odd immersion. Let $\toi(G)$ denote the maximum $t$ such that $G$ contains a totally odd strong immersion of $K_t$. While the class of all graphs $G$ with $\im(G)\le t$ is a sparse class (one with bounded degeneracy), it is not hard to see that if $B$ is a complete bipartite graph, then $\toi(B)\le 2$. Therefore, the following conjecture (which is inspired by one of Churchley \cite{churchley2017odd}) extends Conjecture~\ref{conj:imm} to arbitrarily dense graph classes. 

\begin{conjecture}[Jim\'enez, Quiroz and Thraves Caro~\cite{JimenezQuirozThraves}]\label{conj:church}
    For every graph $G$,  we have $\chi(G)\le \toi(G)$.
\end{conjecture}

Extending the study of totally odd immersions, we consider a question analogous to Question~\ref{qus:im} for the parameter $\toi$. The answer to this question is the following.

\begin{theorem}\label{qus:toi}
    Let $G$ and $H$ be graphs with $\toi(G) = t$ and $\toi(H) = r$, and let $\ast$ be any of the four standard graph products. Then, we have $\toi(G \ast H) \geq \toi(K_t \ast K_r)$.
\end{theorem}

In Sections~\ref{sec:direct} and \ref{sec:cartesian} we show that Theorem \ref{qus:toi} holds for the case of the direct and the cartesian product, respectively. The proof of Theorem 3 in~\cite{collins2023clique} already proves Theorem \ref{qus:toi} for the lexicographic product. Similarly, of Theorem 6 in \cite{wu2025note} gives Theorem~\ref{qus:toi} for the strong product.

We apply Theorem \ref{qus:toi} to give evidence of Conjecture~\ref{conj:church} as follows.  A strengthening of Conjecture~\ref{conj:imm}, different from Conjecture~\ref{conj:church}, had been proposed earlier by Haj\'os \cite{Hajos}, who conjectured that every graph $G$ contains a subdivision of a complete graph on $\chi(G)$ vertices. Haj\'os' conjecture was disproved by Catlin~\cite{C79}, who showed that it fails whenever $\chi(G)\ge 7$. Catlin's counterexamples are remarkably simple, and in particular they can be seen as a particular type of graph product (lexicographic product) between a cycle and a complete graph.  While Haj\'os' conjecture and Conjecture~\ref{conj:church} are incomparable it is interesting to check whether Catlin's counterexamples to Haj\'os' conjecture could be counterexamples to Conjecture~\ref{conj:church}. This possibility was discarded by Jim\'enez, Quiroz and Thraves Caro~\cite{JimenezQuirozThraves}. Yet it is natural to wonder whether graph products can produce minimal counterexamples to Conjecture~\ref{conj:church}. This is, however, discarded by the following.

\begin{corollary}
Let $\ast$ be any of the four standard graph products. No minimum counterexample to Conjecture~\ref{conj:church} is of the form $G \ast H$.
\end{corollary}

We prove now that this is the case for the lexicographic product; the same proof works for the strong product. For the other products, the reader can check the corresponding sections of this paper. As mentioned earlier, by the proof of \cite[Theorem 3]{collins2023clique} we have that  $\toi(G\circ H)\ge \toi(G)\toi(H)$ for every $G$ and $H$. So suppose for a contradiction that there are graphs $G$ and $H$ so that $G\circ H$ is a minimal counterexample to Conjecture~\ref{conj:church}. By minimality we have $\chi(G)\le \toi(G)$ and $\chi(H)\le \toi(H)$. It is well known that $\chi(G\circ H)\le \chi(G)\chi(H)$, and so we obtain, $\chi(G\circ H)\le\toi(G\circ H)$, a contradiction.

 We end this section with the definitions of the graph products which are the object of our study.

All graphs in this paper are finite, simple and loopless. Let $G$ and $H$ be graphs. The \emph{lexicographic product} graph between $G$ and $H$ is the graph on vertices $V(G)\circ V(H)$ and where two vertices $(g_1,h_1)$, $(g_2,h_2)$ are adjacent if $g_1g_2\in E(G)$, or $g_1=g_2$ and $h_1h_2\in E(H)$. The resulting graph is denoted $G\circ H$. The \emph{Cartesian product} between $G$ and $H$, denoted by $G\Square H$, is the graph on vertices $V(G)\times V(H)$ where two vertices $(g_1,h_1),(g_2,h_2)$ are adjacent if $g_1=g_2$ and $h_1h_2\in E(H)$, or $h_1=h_2$ and $g_1g_2\in E(G)$.
The \emph{direct product} (also known as tensor product) between $G$ and $H$, denoted by $G\times H$, is the graph with vertex set $V(G)\times V(H)$ and where two vertices $(g_1,h_1),(g_2,h_2)$ are adjacent if $g_1g_2\in E(G)$ and $h_1h_2\in E(H)$.
The \emph{strong product} between $G$ and $H$, denoted by $G\boxtimes H$, corresponds to the union $(G\Square H) \cup (G\times H)$. That is, $V(G\boxtimes H)=V(G)\times V(H)$ and two vertices $(g_1,h_1),(g_2,h_2)$ are adjacent if $g_1=g_2$ and $h_1h_2\in E(H)$, or $h_1=h_2$ and $g_1g_2\in E(G)$, or $g_1g_2 \in E(G)$ and $h_1h_2\in E(H)$.

\section{Direct Product} \label{sec:direct}

In this section we study the behaviour of the parameter $\toi$ under the operation of direct product.   The direct product of a bipartite graph with any other graph is a bipartite graph. Thus, $\toi(G \times H)\le 2$, if either $G$ or $H$ is bipartite. Thus, we only need to study the case of $\toi(G \times H)$ when neither $G$ nor $H$ is a bipartite graph.

\begin{theorem} \label{theo:GH_direct_KtKs}
Let $G$ and $H$ be two graphs with $\toi(G) = t \geq 3$ and $\toi(H) = s \geq 3$, then
 $$\toi(G \times H) \geq \toi(K_{t} \times K_s). 
  $$ \end{theorem}

 \begin{proof}

Let $u_1,u_2, \dots u_t$ be the terminals of a totally odd strong immersion of $K_t$ in $G$ and, for $1 \leq i<i' \leq t$, let $P_{ii'}$ be the corresponding odd path joining $u_i$ and $u_{i'}$.  Analogously, we let $v_1,v_2, \dots v_s$ be the terminals of a totally odd strong immersion of $K_s$ in $H$ and $Q_{jj'}$ be the corresponding odd path joining $v_j$ and $v_{j'}$, where $1 \leq j<j' \leq s$.  Suppose that $\toi(K_t\times K_s) = r$, and hence, $K_t\times K_s$ has a totally odd immersion of $K_r$. Let $V(K_t)=\{x_1, \ldots, x_t\}$, $V(K_s)=\{y_1, \ldots, y_s\}$, and $c_1, \cdots, c_r \subseteq  \{(x_i,y_j) \mid 1 \leq i \leq t, 1 \leq j \leq s\}$ be the terminals of a 
totally odd immersion of $K_r$ in $K_t\times K_s$.
Let us consider the injective function $\phi: K_t\times K_s \rightarrow G \times H$ as $\phi(x_i,y_j) = (u_i,v_j)$, for all $1 \leq i \leq t$ and $1 \leq j \leq s$. 
Our plan is to construct a totally odd strong immersion of $K_r$ in $G \times H$ with  set of terminals $$\{\phi(c) \in \{(u_i,v_j) \mid 1 \leq i \leq t, 1 \leq j \leq s\} \mid  c \in \{c_1,\ldots, c_r\}\}.$$
For $1 \leq i<i' \leq t$ and $1 \leq j < j' \leq s$, we let $P_{ii'}= u_i-a_1-a_2-\dots-a_{k}-u_{i'}$ and 
 $Q_{jj'}=v_j-b_1-b_2-\dots-b_{\ell}-v_{j'}$, as before.
Since $P_{ii'}$, $Q_{jj'}$ are odd paths,  $k$ and $\ell$ are even numbers. We are going to define two odd paths  $M_{(i,j)-(i',j')}$ and $M_{(i,j')-(i',j)}$ in $G \times H$, connecting $(u_i,v_j)$ to $(u_{i'},v_{j'})$ and $(u_i,v_{j'})$ to $(u_{i'},v_{j})$, respectively, as follows. 

When both $k$ and $\ell$ are zero, the vertices $(u_i,v_j)$ and $(u_{i'},v_{j'})$ are adjacent, and hence  $M_{(i,j)-(i',j')}$ and $M_{(i,j')-(i',j)}$ are defined to be the edges $(u_i,v_j)-(u_{i'},v_{j'})$ and $(u_i,v_{j'})-(u_{i'},v_j')$. Hence, we may assume that at least one of the numbers $k$ or $\ell$ is at least one.
We split the definition of the paths into cases depending on the relation between $k$ and $\ell$.

\noindent\textbf{Case: $0< k < \ell$.} We define $M_{(i,j)-(i',j')}$ through the sequence of vertices

\begin{center}
    $(u_i,v_j)-(a_1,b_1)-(a_2,b_2)-\dots -(a_k,b_k)-(u_{i'},b_{k+1})-\dots -(a_k,b_{\ell})-(u_{i'},v_{j'}),$
\end{center}
where from $(a_k,b_k)$ to  $(u_{i'},v_{j'})$ the first coordinate alternates between $a_k$ and $u_{i'}$, and the second coordinate continues moving along the path $Q_{jj'}$. Note that this is a valid path since $k$ and $\ell$ are of same parity. Now, 
we define $M_{(i,j')-(i',j)}$ through the sequence of vertices
\begin{center}
    $(u_i,v_{j'})-(a_1,b_{\ell})-(a_2,b_{\ell-1})-\dots -(a_{k},b_{\ell-k+1}) - (u_{i'},b_{\ell-k}) -\dots-(a_{k},b_{1})-(u_{i'},v_{j}),$
\end{center}
where from $(a_k,b_{\ell-k+1})$ to  $(u_{i'},v_{j})$ the first coordinate alternates between $a_k$ and $u_{i'}$, and the second coordinate continues  moving along the path $Q_{jj'}$ from $v_{j'}$ to $v_{j}$.

We argue that paths $M_{(i,j)-(i',j')}$ and $M_{(i,j')-(i',j)}$ are edge-disjoint. Actually, the paths are vertex-disjoint. 
In $M_{(i,j)-(i',j')}$ there are five types of vertices: the end vertices $(u_i,v_j)$ and $(u_{i'},v_{j'})$,  vertices $(a_p,b_p)$, vertices $(a_k, b_{even})$, and vertices $(u_{i'}, b_{odd})$. None of these vertices appear in  $M_{(i,j')-(i',j)}$, since 
vertices in $M_{(i,j')-(i',j)}$ are of the following six types: the end vertices $(u_i,v_{j'})$ and $(u_{i'},v_{j})$,  vertices $(a_{odd},b_{even})$, $(a_{even},b_{odd})$, 
vertices $(a_k, b_{odd})$, and vertices $(u_{i'}, b_{even})$.
\vspace{0.3cm}

\noindent\textbf{Case: $0< \ell <k$.}
We define $M_{(i,j)-(i',j')}$ through the sequence of vertices
\begin{center}
    $
 (u_i,v_j)-(a_1,b_1)-(a_2,b_2)-\dots -(a_\ell,b_\ell)-(a_{\ell+1},v_{j'})-\dots (a_k,b_{\ell})-(u_{i'},v_{j'}),
$
\end{center}
where from $(a_\ell,b_\ell)$ to  $(u_{i'},v_{j'})$ the second coordinate alternates between $b_{\ell}$ and $v_{j'}$, and the first coordinate continues moving along the path $P_{ii'}$.
Now, we define $M_{(i,j')-(i',j)}$ through the sequence of vertices 
\begin{center}
    $(u_i,v_{j'})-(a_1,b_{\ell})-(a_2,b_{\ell-1})-\dots -(a_{\ell},b_{1})-(a_{\ell+1},v_{j})-(a_{\ell+2},b_{1})-\dots-(a_{k},b_{1})-(u_{i'},v_{j}),$
\end{center}
where from $(a_\ell,b_1)$ to  $(u_{i'},v_{j})$ the second coordinate alternates between $b_{1}$ and $v_{j}$, and the first coordinate continues moving along the path $P_{ii'}$. 

As in the previous case,
we can argue that paths $M_{(i,j)-(i',j')}$ and $M_{(i,j')-(i',j)}$ are vertex-disjoint. In $M_{(i,j)-(i',j')}$ there are five types of vertices: the end vertices $(u_i,v_j)$ and $(u_{i'},v_{j'})$,  vertices $(a_p,b_p)$, vertices $(a_{even}, b_{\ell})$, and vertices $(a_{odd}, v_{j'})$. None of these vertices appear in  $M_{(i,j')-(i',j)}$, since 
vertices in $M_{(i,j')-(i',j)}$ are of the following six types: the end vertices $(u_i,v_{j'})$ and $(u_{i'},v_{j})$,  vertices $(a_{odd},b_{even})$, $(a_{even},b_{odd})$, 
vertices $(a_{even}, b_{1})$, and vertices $(a_{odd}, v_{j})$.

\vspace{0.3cm}

\noindent\textbf{Cases: either $0 \in \{k, \ell\}$ or $k=\ell$.} Definitions are given according to the rules stated in Table~\ref{tab:paths1.1} and Table~\ref{tab:paths1.2}.

{\renewcommand{\arraystretch}{2}%
\begin{table}[h!]
\centering

\begin{tabular}{|c|l|}
\hline 
\textbf{definition of} $M_{(i,j)-(i',j')}$ & \textbf{condition} \\
\hline
    $(u_i,v_j)-(a_1,b_1)-(a_2,b_2)-\dots -(a_{k},b_{\ell})-(u_{i'},v_{j'})$ &  $k=\ell$\\
    \hline
    $(u_i,v_j) - (u_{i'},b_1)- (u_i,b_2)-\dots - (u_i,b_{\ell})-(u_{i'},v_{j'})$ &  $k=0$, $\ell \geq2$ 
    \\
    \hline
    $(u_i,v_j)-(a_1,v_{j'})-(a_2,v_{j})-\dots -(a_k,v_{j})-(u_{i'},v_{j'})$ & $\ell=0$, $k\geq2$ 
    \\
    \hline
    \end{tabular}
\caption{Rules for the definition of odd paths joining $(u_i,v_j)$ to $(u_{i'},v_{j'})$ in the case that $1 \leq i < i' \leq t$ and $1\leq j<j'\leq s$. In the subcase $k=0, \ell \geq 2$, the first coordinate alternates between $u_{i}$ and $u_{i'}$, and the second coordinate moves along the path $Q_{jj'}$. In the subcase $\ell=0, k \geq 2$, the second coordinate alternates between $v_{j}$ and $v_{j'}$, and the first coordinate moves along the path $P_{ii'}$. }
\label{tab:paths1.1}
\end{table}}

{\renewcommand{\arraystretch}{2}%
\begin{table}[h!]
\centering

\begin{tabular}{|c|l|}
\hline 
\textbf{definition of} $M_{(i,j')-(i',j)}$ & \textbf{condition} \\
\hline
    $(u_i,v_{j'})-(a_1,b_{\ell})-(a_2,b_{\ell-1})-\dots -(a_{k},b_{1})-(u_{i'},v_{j})$ &  $k=\ell$\\
    \hline
    $(u_i,v_{j'})-(u_{i'},b_{\ell})-(u_i,b_{\ell-1})-\dots-(u_{i},v_{j})-(u_{i'},v_j)$ &  $k=0$, $\ell \geq2$  
    \\
    \hline
    $(u_i,v_{j'})- (a_1,v_j)- (a_2,v_{j'})-\dots- (a_k,v_{j'})-(u_{i'},v_{j})$ & $\ell=0$, $k\geq2$ 
    \\
    \hline
    
    \end{tabular}
\caption{Rules for the definition of paths joining $(u_i,v_{j'})$ to $(u_{i'},v_{j})$ in the case that $1 \leq i < i' \leq t$ and $1\leq j<j'\leq s$. In the subcase $k=0, \ell \geq 2$, the first coordinate alternates between $u_{i}$ and $u_{i'}$, and the second coordinate moves along the path $Q_{jj'}$ from $v_{j'}$ to $v_j$. In the subcase $\ell=0, k \geq 2$, the second coordinate alternates between $v_{j'}$ and $v_{j}$, and the first coordinate moves along the path $P_{ii'}$.}
\label{tab:paths1.2}
\end{table}}

Again, in this case, we can argue similarly that for each subcase, paths $M_{(i,j)-(i',j')}$ and $M_{(i,j')-(i',j)}$ are vertex-disjoint. For the subcase $k=\ell$, in $M_{(i,j)-(i',j')}$ there are two types of vertices: the end vertices $(u_i,v_j)$ and $(u_{i'},v_{j'})$  and vertices $(a_p,b_p)$. None of these vertices appear in  $M_{(i,j')-(i',j)}$, since 
vertices in $M_{(i,j')-(i',j)}$ are of the following three types: the end vertices $(u_i,v_{j'})$ and $(u_{i'},v_{j})$, and  vertices $(a_{odd},b_{even})$, $(a_{even},b_{odd})$. For the subcase $k=0, \ell \geq 2$, in $M_{(i,j)-(i',j')}$ there are three types of vertices: the end vertices $(u_i,v_j)$ and $(u_{i'},v_{j'})$, and vertices $(u_{i'}, b_{odd})$, $(u_{i}, b_{even})$. None of these vertices appear in  $M_{(i,j')-(i',j)}$, since 
vertices in $M_{(i,j')-(i',j)}$ are of the following three types: the end vertices $(u_i,v_{j'})$ and $(u_{i'},v_{j})$,  and vertices $(u_{i'}, b_{even})$, $(u_{i}, b_{odd})$. Finally, for the subcase $\ell=0, k \geq 2$, in $M_{(i,j)-(i',j')}$ there are three types of vertices:
the end vertices $(u_i,v_j)$ and $(u_{i'},v_{j'})$, and vertices $(a_{odd}, v_{j'})$, $(a_{even}, v_{j})$. None of these vertices appear in  $M_{(i,j')-(i',j)}$, since 
vertices in $M_{(i,j')-(i',j)}$ are of the following three types: the end vertices $(u_i,v_{j'})$ and $(u_{i'},v_{j})$,  and vertices $(a_{odd}, v_{j})$, $(a_{even}, v_{j'})$.

\vspace{0.3cm}

We claim that the set of all odd paths $M_{(i,j)-(i',j')}$, $M_{(i,j')-(i',j)}$ with $1 \leq i \leq i' \leq t$ and $1\leq j\leq  j'\leq s$ is pairwise edge-disjoint. 
%{\color{red}{Poner argumento limpio de esta afirmación, lo dejé comentado, hay que explicarlo bien.}}

Let us take the path $M_{(i,j)-(i,j')}$ for some $1 \leq i \leq t$ and $1 \leq j < j'\leq  s$. Note that every edge $(x,y)-(x',y')$ in the path $M_{(i,j)-(i,j')}$ satisfies that $x=x'=u_i$ and $yy'\in Q_{jj'}$. Thus, if $M_{(i,j)-(i,j')}$ shares an edge with another path from $\mathcal{M}$, this path should be $M_{(i,q)-(i,q')}$ for some  $1 \leq q < q'\leq s$ with $\{q,q'\} \neq \{j,j'\}$. But then, $yy'$ is an edge in common in $Q_{jj'}$ and $Q_{qq'}$, which contradicts the fact that $Q_{jj'}$ and $Q_{qq'}$ are paths from a totally odd strong immersion of $K_s$ in $H$. Analogously, the same holds for any path $M_{(i,j)-(i',j)}$ for every $1 \leq i<i' \leq t$ and $1 \leq j \leq  s$.

It is left to show that the set of paths

\begin{center}
    $\{M_{(i,j)-(i',j')}, M_{(i,j')-(i',j)}  \mid 1 \leq i < i' \leq t, 1 \leq j < j' \leq s \}$
\end{center}
is pairwise edge disjoint.
Let $M_{(i_1,j_1)-(i_1',j_1')}$ and $M_{(i_2,j_2)-(i_2',j_2')}$ distinct paths with $1 \leq i_1 < i_1' \leq t, 1 \leq i_2 < i_2' \leq t, 1 \leq j_1, j_1'\leq s$ and $1\leq j_2, j_2'\leq s$, with $j_1\neq j_1'$ and $j_2\neq j_2'$.
By the sake of contradiction, we suppose $(x,y)-(x',y')$ is a common edge of the aforementioned paths. Due to the construction of the paths we have that $xx' \in E(P_{i_1i_1'})$ and $yy' \in E(Q_{j_1'j_1})$ 
and at the same time $xx' \in E(P_{i_2i_2'})$ and $yy' \in E(Q_{j_2'j_2})$. However, $P_{i_1 i_1'}$ and $P_{i_2 i_2'}$ ($Q_{j_1 j_1'}$ and $Q_{j_2 j_2'}$, respectively) are paths of the totally odd strong immersion of $K_t$ in $G$ ($K_s$ in $H$, respectively), and hence, they are either edge-disjoint or equal. Thus,  $P_{i_1 i_1'}=P_{i_2 i_2'}$ and  $Q_{j_1 j_1'}=Q_{j_2 j_2'}$, and since $M_{(i_1,j_1)-(i_1',j_1')}$ and $M_{(i_2,j_2)-(i_2',j_2')}$  are distinct, then, $M_{(i_2,j_2)-(i_2',j_2')}$ is the path $M_{(i_1,j_1')-(i_1',j_1)}$. However, we have already proven that  $M_{(i,j)-(i',j')}$ and $M_{(i,j')-(i',j)}$ are edge-disjoint (actually, vertex-disjoint)
for every $1 \leq i < i' \leq t$ and $1 \leq j < j'\leq s$.

Now,  we show that 
$\{\phi(c) \in \{(u_i,v_j) \mid 1 \leq i \leq t, 1 \leq j \leq s\} \mid  c \in \{c_1,\ldots, c_r\}\}$ is the set of terminals of a totally odd strong immersion of $K_r$ in $G \times H$.

Suppose $\phi(c)=(u_i,v_j)$ and $\phi(c')=(u_{i'},v_{j'})$ for some $ 1 \leq i \leq i' \leq t$ and $j, j' \in \{1, \ldots, s\}$.

 \begin{itemize}
     \item If $i \neq i'$ and $j \neq j'$, then for the path joining $\phi(c)$ and $\phi(c')$ we consider
     $M_{(i,j)-(i',j')}$ or $M_{(i,j')-(i',j)}$ depending on whether $j<j'$ or $j>j'$. 
     
     \item If either $i=i'$ or $j=j'$, then first we take the odd path $$(x_i,y_j)-(x_{i_1},y_{j_1})- (x_{i_2},y_{j_2}) -\dots- (x_{i_{q}},y_{j_{q}})-(x_{i'},y_{j'})$$ 
     joining $(x_i,y_j)$ to $(x_{i'},y_{j'})$ in the totally odd strong immersion of $K_r$ in $K_t \times K_s$. Note that each edge 
     $(x,y)- (x',y')$ from this path satisfies that $x\neq x'$ and $y\neq y'$ due to the definition of direct product. Recall that $\phi(x_i,y_j)= (u_i,v_j)$ for $1 \leq i \leq t$ and $1 \leq j \leq s$. See Figure \ref{fig:direct1} for an illustration of this case.
     Thus, we can now define the odd path of the immersion of $K_r$ in $G \times H$
     joining  $(u_i,v_j)$ to $(u_{i'},v_{j'})$ as the union of the following odd paths $$M^*_{(i,j)-({i_1},{j_1})}-M^*_{(i_1,j_1)-(i_2,j_2)}-M^*_{({i_2},{j_2})-({i_3},{j_3})}- \dots 
  -M^*_{(i_{q-1},j_{q-1})-(i_{q},j_{q})}-M^*_{(i_{q},j_{q})-(i',j')}$$
where 
  $$M^*_{(h,g)-(h',g')} = \begin{cases}
    M_{(h,g)-(h',g')} ~&\text{ if } h<h' \text{ and } g<g' \\
    M_{(h',g')-(h,g)} ~&\text{ if } h>h' \text{ and } g'<g
\\
    M_{(h,g')-(h',g)} ~&\text{ if } h<h' \text{ and } g>g'
    \\
    M_{(h',g)-(h,g')} ~&\text{ if } h>h' \text{ and } g'>g'    
  \end{cases}$$
  
 \end{itemize}

Note that these odd paths are edge-disjoint as $$\{M_{(i,j)-(i',j')}, M_{(i,j')-(i',j)}  \mid 1 \leq i < i' \leq t, 1 \leq j < j' \leq s \}$$ is a set of pairwise edge-disjoint odd paths in $G \times H$, as proven in the beginning of the proof. Therefore, $G \times H$ has $K_r$ as a totally odd strong immersion. 

\begin{figure}[h!]
    \centering
    \includegraphics[scale = 0.55]{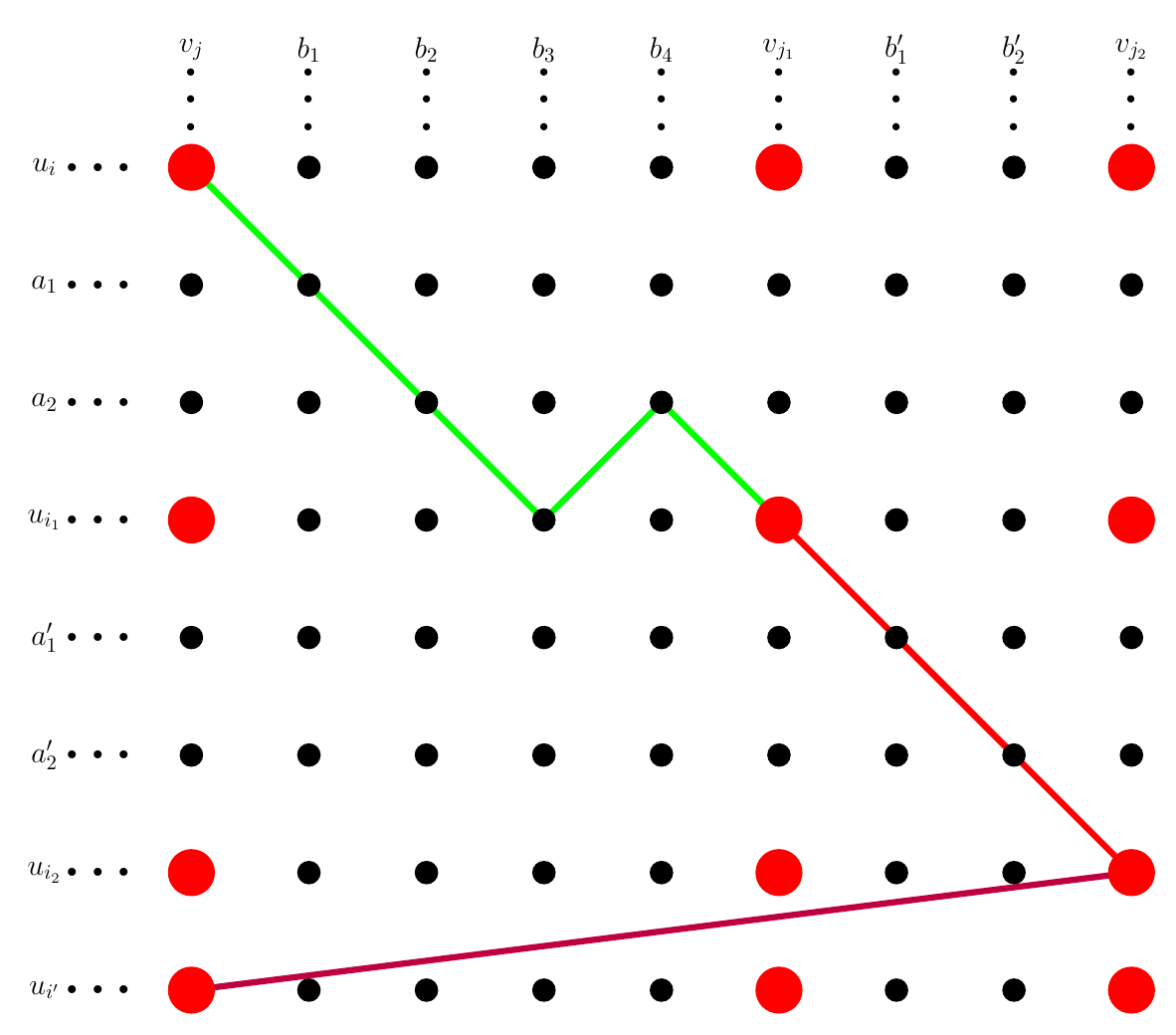}
    \caption{Let us suppose that the odd path of the totally odd strong immersion of $K_r$ in $K_t\times K_s$ that joins $(x_i,y_j)$ and $(x_{i'},y_j)$ is given by the sequence of vertices $(x_i,y_j)-(x_{i_1},y_{j_1})- (x_{i_2},y_{j_2}) -(x_{i'},y_{j})$. The illustration shows the result of the union of the paths $M^*_{(i,j)-({i_1},{j_1})}$ in green, $M^*_{({i_1},{j_1})-({i_2},{j_2})}$ in red and $M^*_{({i_2},{j_2})-(i',j)}$ in dark red, which corresponds to the odd path joining $(u_i,v_j)$ and $(u_{i'},v_j)$ in the totally odd strong immersion of $K_r$ in $G\times H$.}
    \label{fig:direct1}
\end{figure}

\end{proof}

\begin{corollary}
      The direct product of two graphs is not a minimum counterexample to Conjecture \ref{conj:church}.
\end{corollary}
\begin{proof}
 Let us assume that there are graphs $G$ and $H$ so that $G\times H$ is a minimum counterexample to Conjecture \ref{conj:church}, and that $\toi(G)=t, \toi(H)=s$. Since $\times$ is commutative, we can assume that $t \le s$. 
Due to Theorem~\ref{theo:GH_direct_KtKs}, we have that $\toi(G\times H) \geq \toi(K_t\times K_s)$. Note that 
 If $V(K_t) =\{x_1,x_2,x_3,\dots,x_t\} $, and $V(K_s) = \{y_1,y_2,y_3,\dots,y_s\}$, then the set $\{(x_i,y_i) \in V(K_t\times K_s) :i\in \{1,\dots, t\}\}$ induces $K_t$ as a subgraph in $K_t \times K_s$, and hence $\toi(K_t\times K_s) \geq t =\min\{\toi(G),\toi(H)\}$. Since $G$ and $H$ are not counterexamples to Conjecture \ref{conj:church}, we have $\toi(G) \geq \chi(G)$, $\toi(H) \geq \chi(H)$, and hence $\min\{\toi(G),\toi(H)\} \geq \min\{\chi(G),\chi(H)\}$.  Furthermore, we have that $\min\{\chi(G),\chi(H)\} \geq \chi(G\times H)$, which finally yields a contradiction.
\end{proof}

Our next result regarding the paramenter $\toi$ for the direct product of graphs tells us that the largest totally odd strong immersion in $K_{2t} \times K_s$ of a complete graph has at least $ts$ vertices. We think, it is possible to improve the bound of Theorem \ref{teo:directKst}, however many more ad-hock paths would be needed.
In~\cite{collins2023clique}, the authors obtain (in Theorem 18) the largest strong immersion of the product of two complete graphs, where most of the paths are of even length. Here, we construct paths of length one and three.

\begin{theorem}\label{teo:directKst}
    $\toi(K_{2t} \times K_s) \geq t  s$, for all natural numbers $t \geq 6$ and $s \geq 5$.
\end{theorem}

\begin{proof}
Let $\{u_1, u_2, u_3, \dots u_{2t}\}$ and $\{v_1, v_2, v_3\dots v_s\}$ be the set of vertices in $K_{2t}$ and $K_{s}$, respectively. We claim that $K_{2t} \times K_{s}$ has a $K_{ts}$ as a totally odd strong immersion with the set of terminals $A := \{(u_{2i-1},v_{j}) \mid 1 \leq i \leq t, 1 \leq j \leq s\}$. First, note that two elements of $A$ that differ in both coordinates are adjacent in the product graph. Let us denote this set of edges by $\mathcal{E}_s$. Thus, we only need to connect pairs of vertices that differ in exactly one coordinate. This connection is made through edge-disjoint paths of length three, and we use edges in $E(K_{2t} \times K_{s})\setminus \mathcal{E}_s$. In order to prove that the paths that we define are pairwise edge-disjoint, we partition $E(K_{2t} \times K_{s})\setminus \mathcal{E}_s$ into classes according to the function
$f: E(K_{2t} \times K_{s})\setminus \mathcal{E}_s \longrightarrow \mathbb{Z}^{3}$ defined as $$f(u_i,v_j)(u_{i'},v_{j'}) = (i~\mod(2),i-i',j'-j) \quad \text{whenever} \quad j<j'.$$ 

An example is illustrated in Figure \ref{def:tuples}.
\begin{figure}[ht!]
    \centering
\tikzset{every picture/.style={line width=0.75pt}} %set default line width to 0.75pt        

\begin{tikzpicture}[x=0.75pt,y=0.75pt,yscale=-1,xscale=1,scale=2]
%uncomment if require: \path (0,748); %set diagram left start at 0, and has height of 748

%Shape: Circle [id:dp04534480720073408] 
\draw  [fill={rgb, 255:red, 0; green, 0; blue, 0 }  ,fill opacity=1 ] (41.23,27.48) .. controls (41.23,24.79) and (43.42,22.6) .. (46.12,22.6) .. controls (48.81,22.6) and (51,24.79) .. (51,27.48) .. controls (51,30.18) and (48.81,32.37) .. (46.12,32.37) .. controls (43.42,32.37) and (41.23,30.18) .. (41.23,27.48) -- cycle ;
%Shape: Circle [id:dp5314038095616349] 
\draw  [fill={rgb, 255:red, 0; green, 0; blue, 0 }  ,fill opacity=1 ] (81.63,67.48) .. controls (81.63,64.79) and (83.82,62.6) .. (86.52,62.6) .. controls (89.21,62.6) and (91.4,64.79) .. (91.4,67.48) .. controls (91.4,70.18) and (89.21,72.37) .. (86.52,72.37) .. controls (83.82,72.37) and (81.63,70.18) .. (81.63,67.48) -- cycle ;
%Straight Lines [id:da699617646888828] 
\draw    (46.12,27.48) -- (86.52,67.48) ;
%Shape: Circle [id:dp03215768335816338] 
\draw  [fill={rgb, 255:red, 0; green, 0; blue, 0 }  ,fill opacity=1 ] (192.92,26.27) .. controls (192.92,23.57) and (195.1,21.38) .. (197.8,21.38) .. controls (200.5,21.38) and (202.68,23.57) .. (202.68,26.27) .. controls (202.68,28.96) and (200.5,31.15) .. (197.8,31.15) .. controls (195.1,31.15) and (192.92,28.96) .. (192.92,26.27) -- cycle ;
%Shape: Circle [id:dp8350081960127034] 
\draw  [fill={rgb, 255:red, 0; green, 0; blue, 0 }  ,fill opacity=1 ] (152.7,66.42) .. controls (152.7,63.72) and (154.89,61.53) .. (157.58,61.53) .. controls (160.28,61.53) and (162.47,63.72) .. (162.47,66.42) .. controls (162.47,69.11) and (160.28,71.3) .. (157.58,71.3) .. controls (154.89,71.3) and (152.7,69.11) .. (152.7,66.42) -- cycle ;
%Straight Lines [id:da4998126784552368] 
\draw    (197.8,26.27) -- (157.58,66.42) ;
%Straight Lines [id:da5179487684085782] 
\draw  [dash pattern={on 4.5pt off 4.5pt}]  (43.85,81.11) -- (91.85,81.11) ;
%Straight Lines [id:da33787020160008163] 
\draw  [dash pattern={on 4.5pt off 4.5pt}]  (158.92,80.04) -- (200.25,80.04) ;
%Straight Lines [id:da7886523134140213] 
\draw  [dash pattern={on 4.5pt off 4.5pt}]  (101.85,23.11) -- (102.52,71.78) ;
%Straight Lines [id:da5008437781010902] 
\draw  [dash pattern={on 4.5pt off 4.5pt}]  (212.92,22.71) -- (213.59,71.38) ;

% Text Node
\draw (62.25,84.51) node [anchor=north west][inner sep=0.75pt]  [font=,xscale=0.8,yscale=0.8]  {$5$};
% Text Node
\draw (178.12,84.64) node [anchor=north west][inner sep=0.75pt]  [font=,xscale=0.8,yscale=0.8]  {$7$};
% Text Node
\draw (106.25,44.38) node [anchor=north west][inner sep=0.75pt]  [font=,xscale=0.8,yscale=0.8]  {$1$};
% Text Node
\draw (217.59,45.44) node [anchor=north west][inner sep=0.75pt]  [font=,xscale=0.8,yscale=0.8]  {$2$};
% Text Node
\draw (52,20.51) node [anchor=north west][inner sep=0.75pt]  [font=,xscale=0.8,yscale=0.8]  {$(u_{1},v_{5})$};
% Text Node
\draw (52.25,64.11) node [anchor=north west][inner sep=0.75pt]  [font=,xscale=0.8,yscale=0.8]  {$(u_{2},v_{10})$};

% Text Node
\draw (165,20.51) node [anchor=north west][inner sep=0.75pt]  [font=,xscale=0.8,yscale=0.8]  {$(u_{2},v_{8})$};
% Text Node
\draw (164.25,60.11) node [anchor=north west][inner sep=0.75pt]  [font=,xscale=0.8,yscale=0.8]  {$(u_4,v_1)$};
% Text Node
\draw (66.4,38.8) node [anchor=north west][inner sep=0.75pt]  [font=,xscale=0.8,yscale=0.8]  {$e_{1}$};
% Text Node
\draw (170.4,38) node [anchor=north west][inner sep=0.75pt]  [font=,xscale=0.8,yscale=0.8]  {$e_{2}$};
% Text Node
\draw (31.2,102.2) node [anchor=north west][inner sep=0.75pt]  [font=,xscale=0.8,yscale=0.8]  {$f( e_{1}) \ =\ ( 1,-1,5)$};
% Text Node
\draw (151.2,101.8) node [anchor=north west][inner sep=0.75pt]  [font=,xscale=0.8,yscale=0.8]  {$f( e_{2}) \ =\ ( 0,2,7)$};

\end{tikzpicture}
    \caption{The edge $e_1$ in the left hand side belongs to the class $(1,-1,5)$ and the edge $e_2$ in the right hand side belongs to the class $(0,2,7)$.}\label{def:tuples}
\end{figure}
Since $f$ is well-defined, we have that
\begin{claim} \label{claim:trivial}
If $f(e) \neq f(e')$, then $e\neq e'$. In other words, edges in different classes are distinct.
\end{claim}

Furthermore, each non-empty class defined according to the function $f$ contains all translations of an edge. In this context, let $e = (u_{i},v_j)(u_{i'},v_{j'})$ and $e' = (u_{k},v_l)(u_{k'},v_{l'})$ be edges, we say that $e'$ is a \textbf{translation} of $e$ if there exist two integers $a,b\geq0$ such that $k = i +2a$, $k' = i'+2a$, $l = j+b$ and $l' = j'+b$.

\begin{claim} \label{claim:translation}
    Let $e,e'$ be two edges in $E(K_{2t} \times K_{s})\setminus \mathcal{E}_t$. We have that $f(e) = f(e')$ if only if either $e=e'$ or $e'$ is a translation of $e$. 
\end{claim}
\begin{proof}
    The only if direction is trivial. Suppose that $f(e) = f(e')$, if $e=e'$ we are ready. Suppose that $e\neq e'$, we have to prove that $e'= (u_{k},v_l)(u_{k'},v_{l'}) $ is a translation of $e = (u_{i},v_j)(u_{i'},v_{j'})$. First, $i~\mod(2) = k~\mod(2)$, it means that there is a $x\geq0$ such that $i = k+2x$. Now, we have $i'-i=k'-k$, taking $~\mod(2)$ on both sides, we get $i'~\mod(2) = k'~\mod(2)$, then there is a $y$ such that $i' = k'+2y$, it is not hard to see that $x=y$. Suppose that there is no integer $b\geq 0$ such that $j'-j\neq b$ and $l'-l\neq b$, then $j'-j\neq l'-l$, contradiction. Finally, $e'$ is a translation of $e$.
 \end{proof}

%An important observation, and the essential of the next proof is that if we want to prove that two edges with the same value of $f$ are different, then we have to prove that they are translations. On the other hand, if they have different values in $f$, then they are distinct. We only have to care for the edges with the same value in $f$. First we prove that the paths from $\p$ are edge-disjoint, after that we do the same for each $\q_i$, and finally we prove that the paths from $\p$ and $\q_i$ are pairwise edge-disjoint.

In the following, we describe the paths for the totally odd strong immersion. We divide the description into cases. In each case, we define a general "pattern/shape" for the paths (given by the classes of their edges) and take translations of them to cover all pairs of vertices within the case. When in the translation, we move the edges of the paths "too much to the right/left" or "too much down/up", we end up having edges which turn around the grid and as a consequence their classes change; that is the main reason for having sub-cases in each case. However, it is convenient for the reader to keep in mind that the "pattern/shape" of the paths is kept. Furthermore, in each case, we make sure that the edges of the paths described in different cases belong to different classes, and hence the paths described in different cases are naturally edge-disjoint.

%\begin{claim}\label{claim:edges-diff}
%If two edges $e = (u_{i},v_j)%(u_{i'},v_{j'})$ and $e'= (u_{k},v_l)(u_{k'},v_{l'}) $ with $j'> j$ and $l'>l$ belong to the same class according to $f$. Then, $e=e'$ if and only if $i=k$ and $j=l$ (i.e. their first vertex is the same).
%is a translation of . First, $i~mod(2) = k~mod(2)$, it means that there is a $x\geq0$ such that $i = k+2x$. 
%belong to the same class (), then  if and only if their first ordered vertex are the same.
%\end{claim}

    \begin{itemize}
    \item[Case 1] For $1 \leq i \leq t-1$ and $1 \leq j < j' \leq s$, we define the path joining $(2i-1,j)$ and $(2i-1,j')$ as
$$P_{(2i-1,j)-(2i-1,j')} := 
        (u_{2i-1},v_{j})-(u_{2i},v_{j'})-(u_{2i+2},v_{j})-(u_{2i-1},v_{j'})$$

     %   We denote by $\p$ the set obtained by the union of these paths.

We find that the first edge of the path $P_{(2i-1,j)-(2i-1,j')} $ is $ (u_{2i-1},v_{j})-(u_{2i},v_{j'}) $ and hence it is in the class 
{\boldmath{$ (1,-1,j'-j)$}}, 
the second edge $(u_{2i},v_{j'})-(u_{2i+2},v_{j})$
is in the class {\boldmath{$(0,2,j'-j)$}} and the third edge $(u_{2i+2},v_{j})-(u_{2i-1},v_{j'})$ is in the class {\boldmath{$(0,3,j'-j)$}}, where $1 \leq j'-j \leq s-1$.

In the case of $i=t$ and $1 \leq j < j' \leq s$, the path joining $(2t-1,j)$ and $(2t-1,j')$ is defined as $$P_{(2t-1,j)-(2t-1,j')} :=(u_{2t-1},v_{j})-(u_{2t},v_{j'})-(u_{2},v_{j})-(u_{2t-1},v_{j'})$$
and the first edge is in the class {\boldmath{$(1,-1,j'-j)$}}, the second edge in the class {\boldmath{$(0,2-2t,j'-j)$}} and the third edge in the class {\boldmath{$(0, 3-2t,j'-j )$}}, where $1 \leq j'-j\leq s-1$.

Let $\p$ denote the class of all described paths. Now, if we consider two edges $e$, $e'$ in different paths $P, P' \in  \p$ that belong to the same class, 
then they are both the first, second or third edges of the respective path. As the paths are different, they are joining different pairs of vertices, and hence the edges are joining different pairs of vertices as well. For a better visualization of the odd paths, see Figure \ref{fig:pathscase1}.

\begin{figure}[H]
    \centering
    \parbox[t]{0.3\textwidth}{
        \centering
        \includegraphics[scale=0.4]{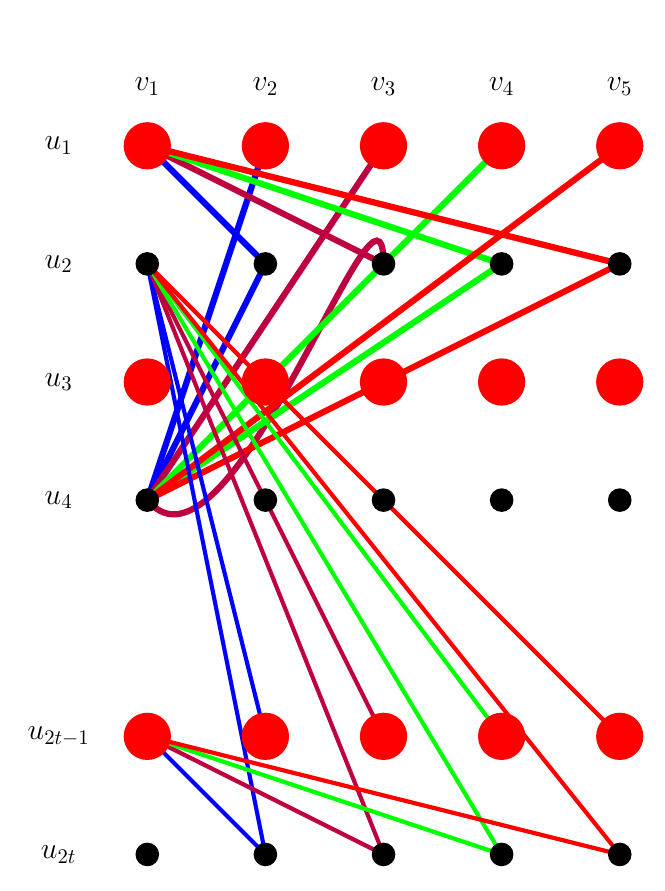}
        \caption{When $j=1$ and $j'\in\{2,3,4,5\}$}
    }
    \hspace{0.02\textwidth}
    \parbox[t]{0.3\textwidth}{
        \centering
        \includegraphics[scale=0.4]{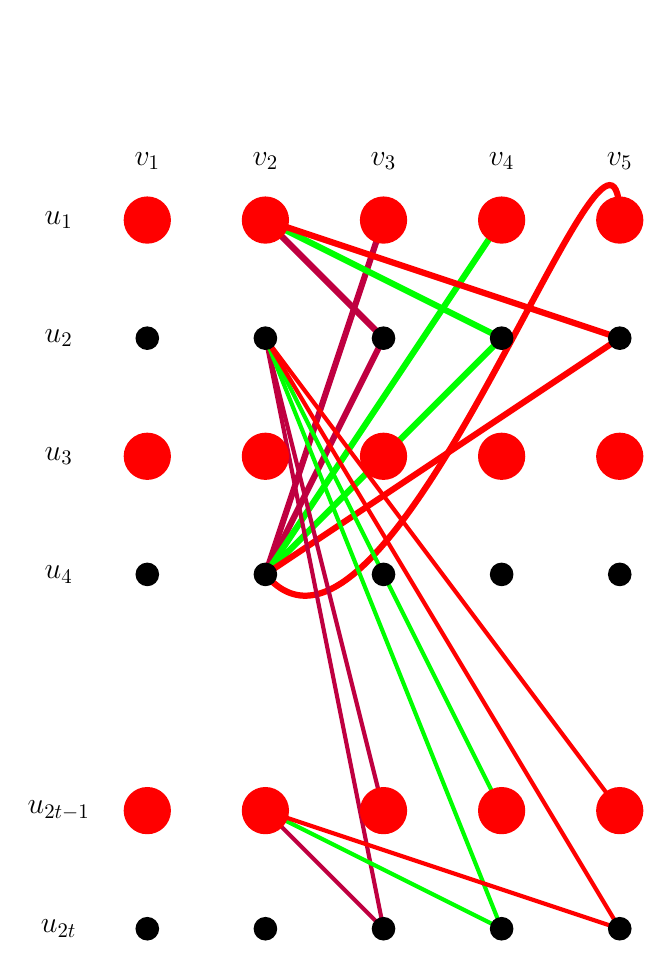}
        \caption{When $j=2$ and $j'\in\{3,4,5\}$}
    }
    \hspace{0.02\textwidth}
    \parbox[t]{0.3\textwidth}{
        \centering
        \includegraphics[scale=0.4]{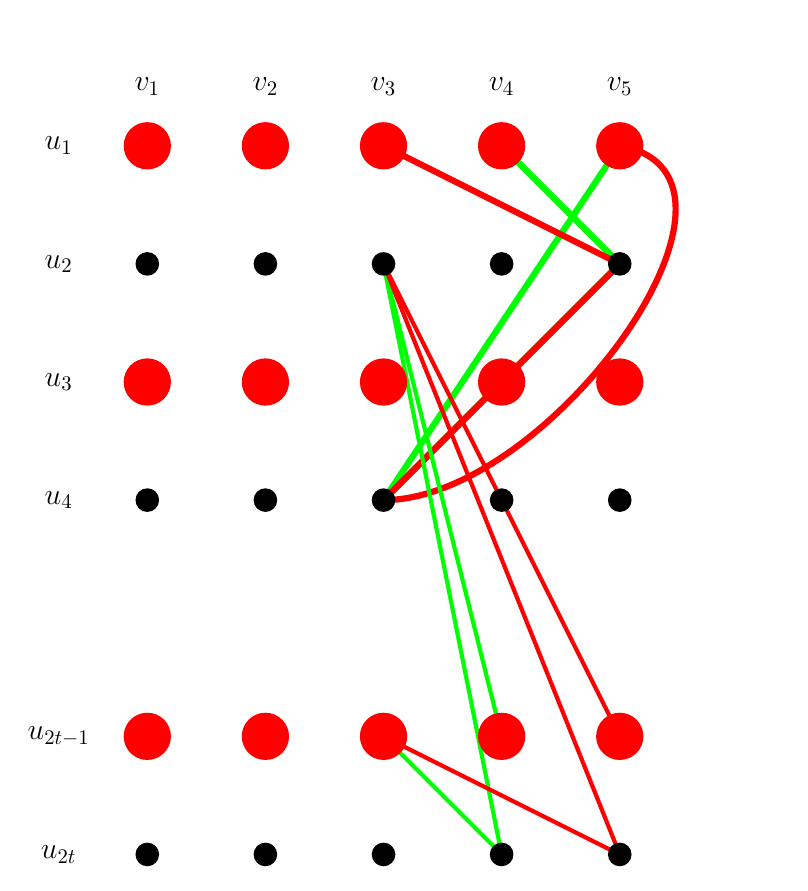}
        \caption{When $j=3$ and $j'\in\{4,5\}$}
    }
    \caption{The picture represents the odd paths in the case 1, when $i\in\{1,t-1\}$, $j, j' \in \{1,2,3,4,5\}$ with $j < j'$. The terminals of the totally odd strong immersion of the $K_{ts}$ are in red while the others are in black.} 
    \label{fig:pathscase1}
\end{figure}

\item[Case 2] For $1\leq i \leq t-1$ and $1 \leq j \leq s-2$, the paths are defined as follows. $$Q_{(2i-1,j)-(2i+1,j)} = (u_{2i-1},v_j)-(u_{2i+2},v_{j+2})-(u_{2i},v_{j+1})-(u_{2i+1},v_j)$$
      
In this case, we define paths so that the edges belong to classes that we have not used yet in the previous case, one can check it exhaustively.
As before, for the described path, we can find the class for each of its edges. The first edge belongs to the class {\boldmath{$(1,-3,2)$}}, the second edges to the class {\boldmath{$(0,-2,1)$}} and the third edge to the class {\boldmath{$(1,1,1)$}}. The class {\boldmath{$(0,-2,1)$}} has not been used in Case 1 because $t>2$.

For $1\leq i \leq t-3$ and $s-1 \leq j \leq s$, the paths are defined as follows.
        $$Q_{(2i-1,j)-(2i+1,j)} = (u_{2i-1},v_j)-(u_{2i},v_{j-1})-(u_{2i+6},v_{j-(s-2)})-(u_{2i+1},v_{j})$$
       % We denote by $\q_4$ the set obtained by the union of these paths.

        The first edges of the paths belong to the class {\boldmath{$(0,1,1)$}}, the second edges of the paths belong to the class {\boldmath{$(0,6,s-3)$}} and the third edge belongs to the class {\boldmath{$(0,5,s-2)$}}.
        Again, the class {\boldmath{$(0,1,1)$}} has not been used in Case 1 because $t>2$.

         For $t-2 \leq i \leq t-1$ and $j \in \{s-1,s\}$, the paths are defined as follows.         $$Q_{(2i-1,j)-(2i+1,j)} = (u_{2i-1},v_j)-(u_{2i-4},v_{j-(s-2)})-(u_{2i},v_{j-(s-3)})-(u_{2i+1},v_{j})$$
        %We denote by $\q_5$ the set obtained by the union of these paths.
        The first edges of the paths belong to the class {\boldmath{$(0,-3,s-2)$}}, the second edges of the paths belong to the class {\boldmath{$(0,-4,1)$}} and the third edge belongs to the class {\boldmath{$(0,-1,s-3)$}}. The three aforementioned classes have not been used in Case 1 because $t>3$.  
        
Again, from the definition it is clear that edges from different paths in the same class are translations of each other (and hence not equal), since they join different pairs of vertices. For an example of some odd path, we refer to the Figure \ref{fig:Paths_in_Q1}.

\begin{figure}[H]
    \centering
\includegraphics[scale=0.5]{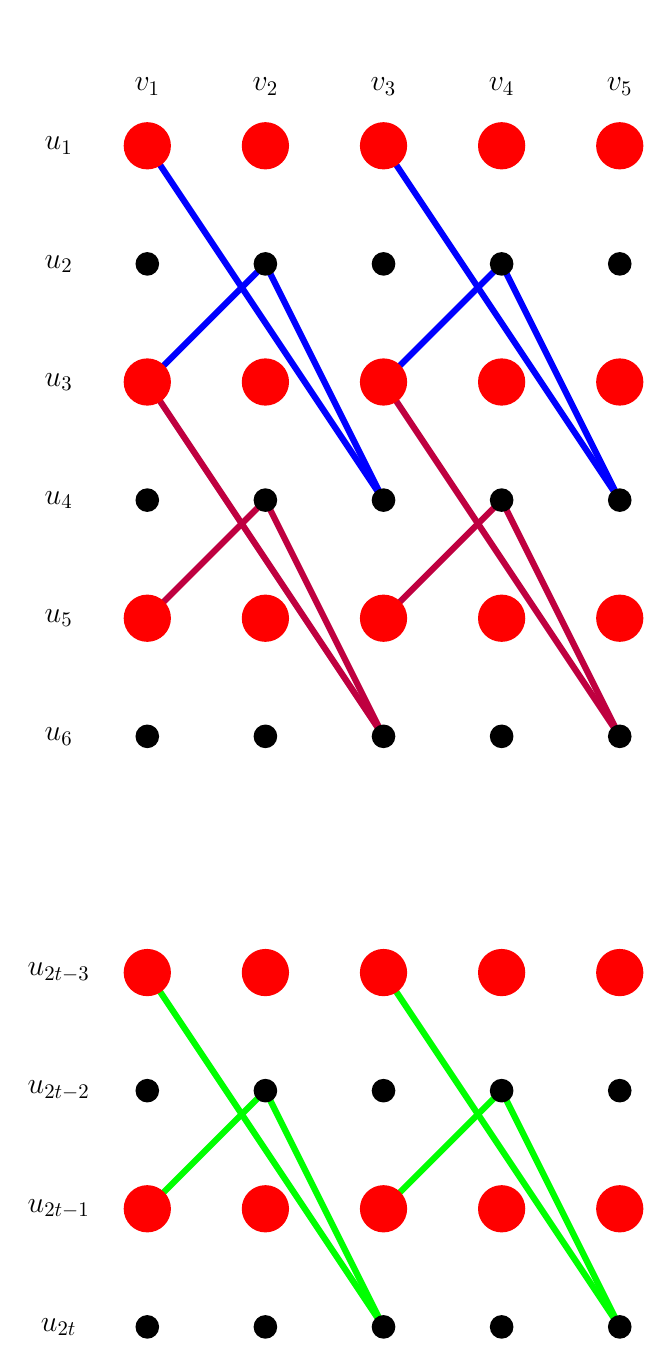}
    \caption{The picture shows the first odd paths described in case 2, when $i\in \{1,2,\dots,t-1\}$ and $j
    \in\{1,3\}$. For instance, the first blue odd path is $Q_{(1,1)-(1,3)} = (u_1,v_1) - (u_4,v_3)-(u_2,v_2)-(u_3,v_1)$.}
    \label{fig:Paths_in_Q1}
\end{figure}

   \item[Case 3] For $1\leq i < i' \leq t $ with $i' \neq  i+1$ and for $1 \leq j \leq s-2$. The paths are defined as follows.
        $$Q_{(2i-1,j)-(2i'-1,j)} = (u_{2i-1},v_j)-(u_{2i'},v_{j+1})-(u_{2i},v_{j+2})-(u_{2i'-1},v_{j})$$

We define paths so that the edges belong to classes that we have not used in cases 1 and 2.  The first edges of the paths belong to classes {\boldmath{$(1,x,1)$}} with $1-2t\leq x \leq -5$ an odd number. The second edge belongs to classes in {\boldmath{$(0,y,1)$}} with $4 \leq y \leq  2t-2$ an even number. The third edge belongs to the class {\boldmath{$(1,z,2)$}} with $5\leq z\leq 2t-3$ an odd number. Due to the parity of $x,y,z$, the definition of the intervals to which $x, y, z$ belong, and the condition that $s\geq 5$, the classes described here have not been considered in previous cases.  %Edges of different paths from $\mathcal{Q}_2$.

    For $1\leq i < i' \leq t $ with $i' \neq  i+1$  and for $j=s-1$, , except for when $i=1$ and $i'=t$ simultaneously. The  paths are as follows.
$$Q_{(2i-1,s-1)-(2i'-1,s-1)} = (u_{2i-1},v_{s-1})-(u_{2i'},v_s)-(u_{2i},v_{1})-(u_{2i'-1},v_{s-1})$$

The first edges of the paths belong to classes {\boldmath{$(1,x,1)$}} with $3-2t\leq x \leq -5$ an odd number. The second edges belong to classes {\boldmath{$(0,y,s-1)$}} with $4-2t\leq y \leq -4$ an even number. The third edges belong to classes {\boldmath{$(0,z,s-2)$}} with $5-2t \leq z\leq -5$ an odd number. Due to the parity of $x,y,z$, the definition of the intervals to which $x, y, z$ belong, and the condition that $s\geq 5$, the classes described here have not been considered in previous cases. However, note that we need a different solution to join the pair $(1,s-1)-(2t-1,s-1)$ since the path $P_{(2t-1,1)-(2t-1,s)}$ defined in Case 1 already uses the edge $(u_{2t},v_{s})-(u_{2},v_{1})$.

 For $1\leq i < i' \leq t $ with $i' \neq  i+1$  and for $j=s$, except for when $i=1$ and $i'=t$ simultaneously. The path is defined as follows.
$$Q_{(2i-1,s)-(2i'-1,s)} = (u_{2i-1},v_{s})-(u_{2i'},v_{1})-(u_{2i},v_{2})-(u_{2i'-1},v_{s})$$

The first edges of the paths belong to classes {\boldmath{$(0,x,s-1)$}} with $5\leq x \leq 2t-3$ an odd number. The second edges belong to classes {\boldmath{$(0,y,1)$}} with $4\leq y \leq 2t-4$  an even number. The third edges belong to the class {\boldmath{$(0,z,s-2)$}} with $5-2t \leq z\leq -5$ an odd number.
As before, because of the parity of $x,y,z$, the definition of the intervals to which $x, y, z$ belong, and the condition that $s\geq 5$, the classes described here have not been considered in previous cases.
Note again that we need a different solution to join the pair $(1,s)-(2t-1,s)$ since the path $P_{(2t-1,2)-(2t-1,s)}$ defined in Case 1 already uses the edge $(u_{2},v_{2})-(u_{2t-1},v_{s})$.

Finally, for $j\in \{s-1,s\}$ and when $i=1$ and $i'=t$ simultaneously. The paths are defined as follows.
        $$Q_{(1,j)-(2t-1,j)} = (u_{1},v_{j})-(u_{2t},v_{j-1})-(u_{2t-6},v_{j-2})-(u_{2t-1},v_{j})$$

       The first edges of the paths belong to the class {\boldmath{$(0,2t-1,1)$}}. The second edges belong to the class in {\boldmath{$(0,-6,1)$}} and the third edges belong to the class {\boldmath{$(0,-5,2)$}}.
   \end{itemize}
\end{proof}

\section{Cartesian Product}\label{sec:cartesian}

In this section, we consider the behaviour of the parameter $\toi$ under the Cartesian product of graphs. Since having $\toi (G)\le 2$ is equivalent to $G$ being bipartite, this means that if $\toi(G)=\toi(H)=2$, then  $\toi(G \Square H)=2=\toi(K_2\Square K_2)$.

So, we focus our attention in the case when at least one of the two graphs is not bipartite, and hence, its $\toi$ value is at least $3$. 
An obvious upper bound for $\toi(K_t \Square K_s)$ is the maximum degree of this product graph $K_t \Square K_s$ plus 1. Note that the degree of every vertex in $K_t \Square K_s$ is $t+s-2$, and hence $\toi (K_t \Square K_s) \le t+s-1$. 
In the next theorem, we show that $t+s-1$ is also a lower bound for $\toi$ of the Cartesian product of two graphs $G$ and $H$ whenever one of these has $\toi$ at least 3.

\begin{theorem}\label{thm:tslargeCart}
    Let $G$ and $H$ be two graphs with  $\toi(G)=t$ and $\toi(H)=s$. If $s \geq 4$, then $\toi(G \Square H) \geq t+s-1$.
\end{theorem}

\begin{figure}[H]
  \centering
  % Crear una estructura de centrado más robusta
  \begin{minipage}[t]{0.4\textwidth}
    \centering
    \includegraphics[scale=1]{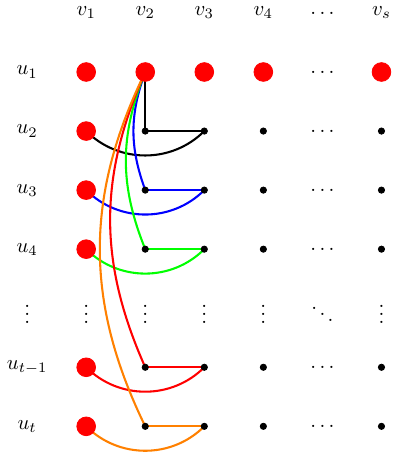}
    \caption{When $i=1$ and $j\in \{1,2,\dots s\}$}
  \end{minipage}
  \hspace{0.1\textwidth} % Espacio controlado
  \begin{minipage}[t]{0.4\textwidth}
    \centering
    \includegraphics[scale=1]{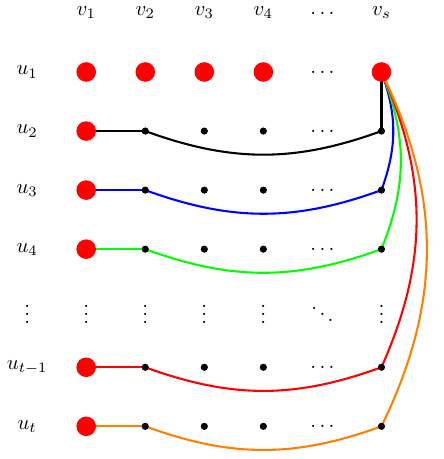}
    \caption{When $i=t$ and $j\in \{1,2,\dots s\}$}
  \end{minipage}
  
  \vspace{1em} % Espacio vertical antes del título general
  \caption{The picture shows the odd paths of the totally odd strong immersion of $K_{t+s-1}$ where the red vertices are its terminals, In particular the odd paths which are described are when $i\in \{2,t\}$ and $j\in \{2,3,\dots,s\}$.}
  \label{fig:Cartesian_Ks_Kt}
\end{figure}

\begin{proof}
Let $u_1,u_2, \dots u_t$ be the terminals of a totally odd strong immersion of $K_t$ in $G$ and, for $1 \leq i<i' \leq t$, let $P_{ii'}$ be the corresponding odd path joining $u_i$ and $u_{i'}$.  Analogously, we let $v_1,v_2, \dots v_s$ be the terminals of a totally odd strong immersion of $K_s$ in $H$ and $Q_{jj'}$ be the corresponding odd path joining $v_j$ and $v_{j'}$, where $1 \leq j<j' \leq s$.

We claim that $G \Square H$ has $K_{s+t-1}$ as a totally odd strong immersion with terminals $T_1 \cup T_2$, where $$T_1 = \{(u_1,v_1),(u_1,v_2),(u_1,v_3),\dots, (u_1,v_s)\}$$ and $$T_2 = \{(u_2,v_1),(u_3,v_1),(u_4,v_1)\dots (u_t,v_1)\}.$$

In the following, we define odd paths 
$M_{(i,j)-(i',j')}$  in $G \Square H$ joining the terminals $(u_i,v_j)$ and $(u_{i'},v_{j'})$, and  show that these paths are pairwise edge-disjoint.
We know that $T_1$ induces 
 a copy of $H$ in the graph $G \Square H$. So, to connect $(u_1,v_i)$ to $(u_1,v_{i'})$, we can use the totally odd strong immersion path in that copy of $H$. More specifically, for any $1 \leq j < j' \leq s$, if $Q_{jj'}$ is $v_j-b_1-b_2-\dots b_{2\ell}-v_{j'}$, then we consider the odd path $$M_{(1,j)-(1,j')}:=(u_1,v_j)-(u_1,b_1)-(u_1,b_2) \dots (u_1,b_{2\ell})-(u_1,v_{j'}).$$ These paths are pairwise edge-disjoint, since $Q_{jj'}$s' are pairwise edge-disjoint. We use a similar idea to connect the terminals in $T_2$. For any $2 \leq i<i' \leq t$, if $P_{ii'}$ is $u_i-a_1-a_2- \dots a_{2k}-u_j$, then we consider the odd path $$M_{(i,1)-(i',1)} := (u_i,v_1)-(a_1,v_1)-(a_2,v_1) \dots (a_{2k},v_1)-(u_{i'},v_1).$$ Again, $P_{ii'}$s' are pairwise edge-disjoint and hence so are $M_{(i,1)-(i',1)}$. Furthermore, any vertex in the path $M_{(1,j)-(1,j')}$ has $u_1$ in the first coordinate. However, no vertex in $M_{(i,1)-(i',1)}$ has $u_1$ in the first coordinate. 

Now we need to connect the vertices of $T_1$ to the vertices of $T_2$. For that, we use the following odd paths $$P_{1j} = u_1-x_1-x_2-\dots -x_{2k}-u_j, \,\, Q_{i(i+1)} = v_i-y_1-y_2-\dots -y_{2k'}-v_{i+1}$$ and $$Q_{(i+1)1} = v_{i+1}-z_1-z_2-\dots -z_{2k''}-v_{1}.$$ For $ 2 \leq i \leq s$ and $2 \leq j \leq t-1 $ , define the odd path $M_{(1,i)-(j,1)}$ as the path obtained by taking the union of the following three paths, 

$$(u_1,v_i)-(x_1,v_i)-(x_2,v_i)- \dots -(x_{2k},v_i)-(u_j,v_i),$$ $$(u_j,v_i)-(u_j,y_1)-(u_j,y_2)- \dots -(u_j,y_{2k'})-(u_j,v_{i+1})$$ and $$(u_j,v_{i+1})-(u_j,z_1)-(u_j,z_2)- \dots -(u_j,z_{2k''})-(u_j,v_1).$$%the path  $P_{(v_1,u_i),(v_j,u_i)}\cup P_{(v_j,u_i),(v_j,u_{i+1})}\cup P_{(v_j,u_{i+1}),(v_j,u_1)}.$

Finally define the odd path $M_{(1,j)-(t,1)}$ as the path obtained by taking the union of the following three paths,

$$(u_1,v_t)-(x_1,v_t)-(x_2,v_t)- \dots -(x_{2k},v_t)-(u_j,v_t),$$

$$(u_j,v_t)-(u_j,y_{2k'})-(u_j,y_{2k'-1})- \dots -(u_j,y_{1})-(u_j,v_{2})$$
and
$$(u_j,v_{2})-(u_j,z_{2k''})-(u_j,z_{2k''-1})- \dots -(u_j,z_{1})-(u_j,v_1).$$

Since all the odd immersion paths ($P_{ii'}$'s and $Q_{jj'}$s') are pairwise edge-disjoint in the respective graphs $G$ and $H$, the new defined paths $M_{(i,j)-(i',j')}$s' are also pairwise edge-disjoint. See the figure \ref{fig:Cartesian_Ks_Kt} for an example of the paths in the Cartesian product between two complete graphs.
\end{proof}

Now, we study the case when $\toi(G),\toi(H)\le 3$. We start with the case when both graphs satisfy this with equality.

\begin{theorem}\label{thm:t3s3Cart}
    If $G$ and $H$ are two connected graphs with $\toi(G)=\toi(H)=3$, then $\toi(G \Square H) \geq 4$. Furthermore, equality holds when $G=H=K_3$.
\end{theorem}

\begin{proof}
    Let $u_1,u_2, u_3$ be the terminals of a totally odd strong immersion of $K_3$ in $G$ and, for $1 \leq i<i' \leq 3$, let $P_{ii'}$ be the corresponding odd path joining $u_i$ and $u_{i'}$.  Analogously, we let $v_1,v_2,v_3$ be the terminals of a totally odd strong immersion of $K_3$ in $H$ and $Q_{jj'}$ be the corresponding odd path joining $v_j$ and $v_{j'}$, where $1 \leq j<j' \leq 3$.
    
    We claim that $G \Square H$ has a $K_4$ as a totally odd strong immersion. To prove this, we consider the set of terminals $\{(u_1,v_1),(u_2,v_1),(u_3,v_1),(u_1,v_2)\}$
 and define pairwise edge-disjoint odd paths 
$M_{(i,j)-(i',j')}$  in $G \Square H$ joining the terminals $(u_i,v_j)$ and $(u_{i'},v_{j'})$.
    For any $1 \leq i < i' \leq 3$, if $P_{ii'}$ is $u_i-a_1-a_2-\dots a_{2k}-u_{i'}$, then to connect $(u_i,v_1)$ to $(u_{i'},v_1)$ we take the path 
    $$M_{(i,1)-(i',1)} := (u_i,v_1)-(a_1,v_1)-(a_2,v_1)-\dots -(a_{2k},v_1)-(u_{i'},v_1).$$
    Note that $P_{12},P_{23}$ and $P_{13}$ are pairwise edge-disjoint and hence so are $M_{(1,1)-(2,1)}, M_{(2,1)-(3,1)}$ and $M_{(1,1),(3,1)}$. Now we need to connect each of them to $(u_1,v_2)$. When $Q_{12}$ is $v_1-b_1-b_2-\dots -b_{2\ell}-u_2$, define $$M_{(1,1)-(1,2)}:=(u_1,v_1)-(u_1,b_1)-(u_1,b_2)- \dots -(u_1,b_{2\ell})-(u_1,v_2).$$ Note that each vertex in $M_{(i,1)-(i',1)}$ has $v_1$ in the second coordinate, for all $1 \leq i < i' \leq 3$, where as no internal vertex in $M_{(1,1)-(1,2)}$ has that property. Thus, $M_{(1,1)-(1,2)}$ is edge-disjoint with $M_{(i,1)-(i',1)}$, for all $1 \leq i < i' \leq 3$.
    
    To connect $(u_i,v_1)$ to $(u_1,v_2)$ we use the paths $P_{1,i} = u_1-x_1-x_2-\dots- x_{2k}-u_i, Q_{2,3} = v_2-y_1-y_2 -\dots- y_{2\ell}-v_3$ and $Q_{3,1} = v_3-z_1-z_2 -\dots -z_{2\ell'}-v_1$, where $2 \leq i \leq 3$. Define $M_{(1,2)-(i,1)}$ as the path obtained by taking the union of the following three paths $$A_i := (u_1,v_2)-(x_1,v_2)-(x_2,v_2)-\dots -(x_{2k},v_2)-(u_i,v_2),$$ $$B_i := (u_i,v_2)-(u_i,y_1)-(u_i,y_2)- \dots -(u_i,y_{2\ell})-(u_i,v_3)$$ and $$C_i := (u_i,v_3)-(u_i,z_1)-(u_i,z_2)-\dots-(u_i,z_{2\ell'})-(u_i,v_2).$$ Note that the second coordinate of $A_i$ is $v_i$, so $A_2$ and $A_3$ are pairwise edge-disjoint. Furthermore, because of the same reason, $A_2$ and $A_3$ are pairwise edge-disjoint with $M_{(1,1)-(2,1)},M_{(1,1)-(3,1)}$ and $M_{(2,1)-(3,1)}$. Moreover, no internal vertex of the following paths $M_{(1,1)-(1,2)}, B_2, B_3, C_1$ and $C_3$ has $v_i$ in the second coordinate. Thus, $M_{(1,1)-(1,2)}, B_2, B_3, C_2$ and $C_3$ are edge-disjoint with $A_2$ and $A_3$. By a similar argument, we can ensure $B_2, B_3, C_2$ and $C_3$ are edge-disjoint with the considered paths. This proves the first part of the statement.

For the second part of the theorem, we need to show that the Cartesian product of two complete graphs of size three, say with vertex sets $\{x_1,x_2,x_3\}$ and $\{y_1,y_2,y_3\}$, does not contain $K_5$ as a totally odd strong immersion. Suppose for the contradiction that there exists a $K_5$ as a totally odd strong immersion. We claim that either there exists  $1 \leq i \leq 3$ such that $(x_i,y_j)$ is a terminal for all $1 \leq j \leq 3$ or there exists  $1 \leq j \leq 3$ such that  $(x_i,y_j)$ is a terminal for all $1 \leq i \leq 3$.

For the sake of contradiction, we assume that there are no $1 \leq i \leq 3$ such that $(x_i,y_j)$ is a terminal for all $1 \leq j \leq 3$ and no $1 \leq j \leq 3$ such that  $(x_i,y_j)$ is a terminal for all $1 \leq i \leq 3$. This is only possible when $\{(x_k,y_\ell),(x_{k},y_{\ell+1}), (x_{k+1},y_\ell),(x_{k+1},y_{\ell+1}), (x_{k+2},y_{\ell+2})\}$ is the set of terminals, for some $k, \ell \in \{1,2,3\}$, where the indices are in module $3$ ($2+2, 3+1,3+2$ are $1,1$ and $2$, respectively).

Without loss of generality we may assume that $k = 1 = \ell$. We have two odd paths (avoiding the terminals as internal vertices) from $(x_3,y_3)$ to $(x_1,y_1)$: $(x_3,y_3)-(x_2,y_3)-(x_1,y_3)-(x_1,y_1)$ and $(x_3,y_3)-(x_3,y_2)-(x_3,y_1)-(x_1,y_1)$. Without loss of generality, we can consider the path $(x_3,y_3)-(x_2,y_3)-(x_1,y_3)-(x_1,y_1)$. Now, to connect $(x_3,y_3)$ to $(x_1,y_2)$, we have only one option $(x_3,y_3)-(x_3,y_1)-(x_3,y_1)-(x_2,y_1)$. Now, we cannot connect $(x_2,y_2)$ with $(x_3,y_3)$. Hence, the claim we made is true: either there exists an $1 \leq i \leq 3$ such that $(x_i,y_j)$ is a terminal for all $1 \leq j \leq 3$ or there exists $1 \leq j \leq 3$  such that $(x_i,y_j)$ is a terminal for all $1 \leq i \leq 3$. Without loss of generality, we may assume that $(x_1,y_1),(x_2,y_1)$ and $(x_3,y_1)$ are three terminals. Again, without loss of generality, we may assume that $(x_1,y_2)$ is the fourth terminal. Note that at least three edges are required to connect $(x_i,y_j)$ with $(x_{i'},y_{j'})$ by an odd path, whenever $i, i'(\neq i) \in \{1,2,3\}$ and $j, j'(\neq j) \in \{1,2,3\}$. Hence at least 10 edges are already used and we have at most 8 more edges left. Thus $(x_2,y_3)$ and $(x_3,y_3)$ are non-terminal vertices. In other words, the fifth terminal is a vertex from $\{(x_1,y_3), (x_2,y_2), (x_3,y_2)\}$. In the last paragraphs, we argue that in any case, we have shortage of edges.

Suppose $(x_1,y_3)$ is the fifth terminal. Notice that we have only two possible candidates to be odd paths from $(x_1,y_2)$ to $(x_2,y_1)$: $(x_1,y_2)-(x_2,y_2)-(x_2,y_3)-(x_2,y_1)$ and $(x_1,y_2)-(x_3,y_2)-(x_2,y_2)-(x_2,y_1)$. If we consider $(x_1,y_2)-(x_2,y_2)-(x_2,y_3)-(x_2,y_1)$, then we have only one option to connect $(x_1,y_2)$ to $(x_3,y_1)$, namely the odd path $(x_1,y_2)-(x_3,y_2)-(x_3,y_3)-(x_3,y_1)$, these paths are depicted in Figure \ref{fig:FirstCase}. Then we cannot connect $(x_1,y_3)$ with $(x_2,y_1)$ and $(x_3,y_1)$.

\begin{figure}[H]
    \centering
    \includegraphics[scale=0.7]{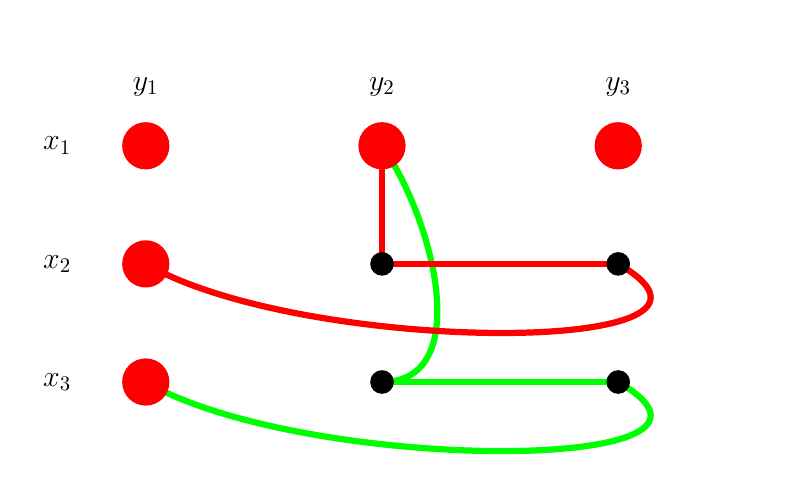}
    \caption{Case when $(x_1,y_3)$ is the fifth terminal. The red path is $(x_1,y_2)-(x_2,y_2)-(x_2,y_3)-(x_2,y_1)$ and the green path 
    is $(x_1,y_2)-(x_3,y_2)-(x_3,y_3)-(x_3,y_1)$.}
    \label{fig:FirstCase}
\end{figure}

Therefore we may assume that either $(x_2,y_2)$ or $ (x_3,y_2)$ is the fifth terminal. Without loss of generality, we may assume that $(x_2,y_2)$ is the fifth terminal. Now to connect $(x_1,y_2)$ with $(x_2,y_1)$ we have only one option $(x_1,y_2)-(x_1,y_3)-(x_2,y_3)-(x_2,y_1)$. Furthermore, we have only one  option to connect $(x_1,y_2)$ to $(x_3,y_1)$, namely $(x_1,y_2)-(x_3,y_2)-(x_3,y_3)-(x_3,y_1)$. Now we cannot connect $(x_2,y_2)$ to $(x_1,y_1)$.
\end{proof}

\begin{figure}[H]
    \centering
   \includegraphics[scale=0.7]{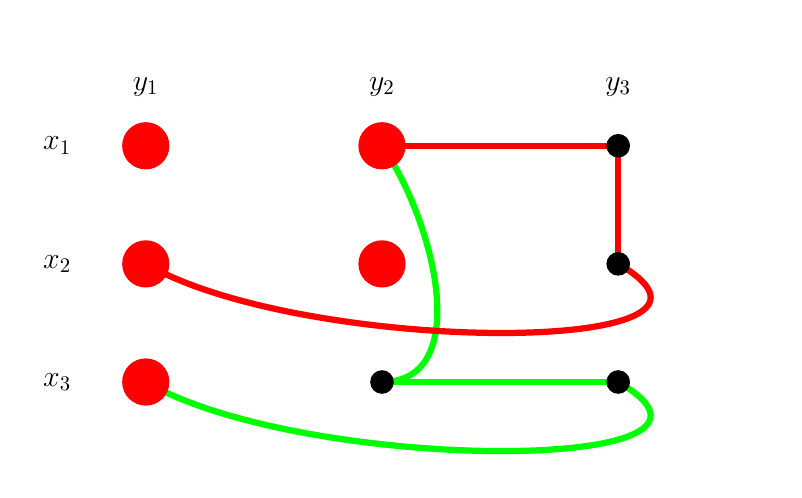}
   \caption{Case when $(x_2,y_2)$ is the fifth terminal. The red path is $(x_1,y_2)-(x_1,y_3)-(x_2,y_3)-(x_2,y_1)$ and the green path 
    is $(x_1,y_2)-(x_3,y_2)-(x_3,y_3)-(x_3,y_1)$.}
   \label{fig:SecondCase}
\end{figure}

\begin{theorem}\label{thm:t3s2Cart}
    Let $G$ and $H$ be two connected graphs with $\toi(G)=3$ and $\toi(H)=2$. 
    Then $\toi(G \Square H) \geq 3$. Furthermore, equality holds when $G=K_3$ and $H=K_2$.
    In general, if $H$ is not isomorphic to $K_2$, then $\toi(G \Square H) \geq 4$.
\end{theorem}

\begin{figure}[H]
    \centering
    \includegraphics[scale=0.7]{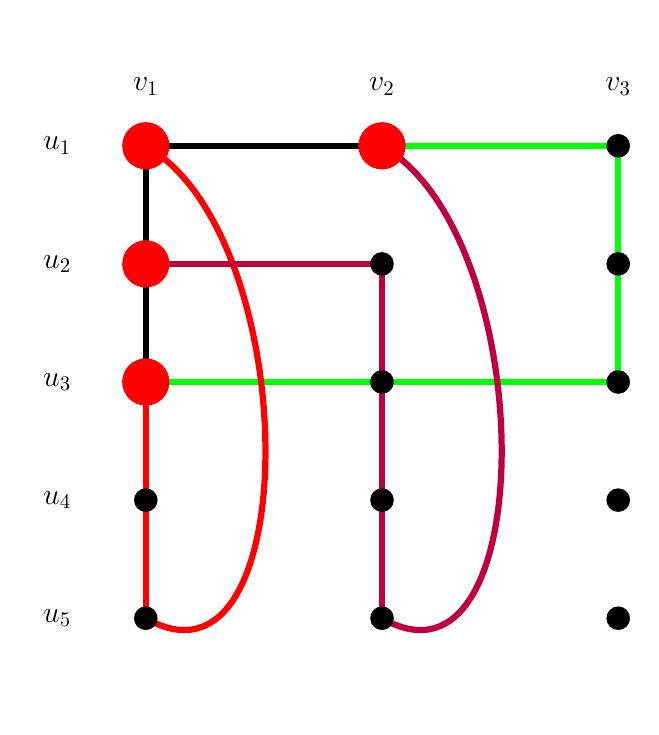}
    \caption{A totally odd strong $K_4$ immersion in $C_5\Square P_3$. Red vertices are the terminals considered in the proof. Black edges are the connections that naturally exist in $C_5\Square P_3$. Red edges represent path $M_{(3,1)-(1,1)}$, while dark red edges represent path $M_{(2,1)-(1,2)}$. Finally, path $M_{(1,3)-(1,2)}$ is depicted in green.}
    \label{fig:odd_imm_k3k2}
\end{figure}

\begin{proof}
Since the $\toi(G)=3$, it is not a bipartite graph. Therefore, $G \Square H$ is a non-bipartite graph and thus, $\toi(G\Square H)$ is at least $3$. Further, $\toi(K_3 \Square K_2)=3$. This finishes the first part of the theorem. For the second part, we assume that $H$ is not $K_2$. Since $H$ is a connected graph, $H$ has a vertex of degree at least two. In other words, $H$ contains a subgraph isomorphic to $P_3$. Again $G$ is not a bipartite graph. Thus, $G$ contains an odd cycle. We need to show that $G \Square H$ has a $K_4$ as a totally odd strong immersion. Because of the previous discussion, it is sufficient to show that $C_{2\ell+1} \Square P_3$ has a $K_4$ as a totally odd strong immersion, for all natural number $\ell$.

Let  $C=C_{2\ell+1}:=u_1-u_2-\dots -u_{2\ell+1}-u_1$ be an odd-cycle and $P_3:=v_1-v_2-v_3$. 
%We consider  can take $\{v_1,v_2,v_3\}$ as the terminals in $C$. 
We claim that $\{(u_1,v_1),(u_2,v_1),(u_3,v_1),(u_1,v_2)\}$ is the set of terminals of a totally odd strong immersion of  $K_4$ in $C_{2\ell+1} \Square P_3$.
In the following, we define pairwise edge-disjoint odd paths $M_{(i,j)-(i',j')}$  in $G \Square H$ joining the terminals $(u_i,v_j)$ and $(u_{i'},v_{j'})$. We define paths the odd path $M_{(1,1)-(1,2)}$,
$M_{(1,1)-(2,1)}$ and $M_{(2,1)-(3,1)}$ to be edges $(u_1,v_1)-(u_2,v_1), (u_1,v_1)-(u_1,v_2)$ and $(u_2,v_1)-(u_3,v_1)$ in $C_{2\ell+1} \Square P_3$ 
(see Figure \ref{fig:odd_imm_k3k2} for an illustration of $C_{2\ell+1} \Square P_3$).

Further, the graph induced in $G \Square H $ by the vertices in $V(C)\times \{v_1\}$ contains an odd path between $(u_3,v_1)$ to $(v_1,u_1)$, namely $$M_{(3,1)-(1,1)}=(u_3,v_1)-(u_4,v_1)-\dots-(u_{2\ell+1},v_1)-(u_1,v_1).$$  %($(v_3,u_1)-(v_4,u_1)- \dots -(v_{2\ell+1},u_1)$ if $\ell >1$, otherwise $(v_3,u_1)$ is adjacent to $(v_1,u_1)$).
To connect $(u_1,v_2)$ to $(u_2,v_1)$  and $(u_1,v_2)$ to  $(u_3,v_1)$ we use the other two copies of $C$ associated with $v_2$ and $v_3$, respectively. The graph induced in $G\Square H$ by the vertices $V(C)\times \{v_2\}$ contains an even path $$(u_2,v_2)-(u_3,v_2)- \dots -(u_{2\ell+1},v_2)-(u_1,v_2)$$ joining $(u_2,v_2)$ to $(u_1,v_2)$. So, we define $M_{(2,1)-(1,2)}$ as
$$(u_2,v_1)-(u_2,v_2)-(u_3,v_2)- \dots -(u_{2\ell+1},v_2)-(u_1,v_2)$$
the odd path joining $(u_2,v_1)$ to $(u_1,v_2)$.

Now, we use the copy of $C$ associated with $v_3$ to get an odd path from $(u_1,v_2)$ to  $(u_3,v_1)$. %Note that the graph induced in $G\Square H$ by the vertices $V(C)\times \{v_3\}$ contains the even path $$(u_1,v_3)-(u_2,v_3)- (u_{3},v_3)$$ joining $(u_2,v_3)$ to $(u_3,v_3)$. 
We define $M_{(1,2)-(3,1)}$ as
$$(u_1,v_2)-(u_1,v_3)-(u_2,v_3)-(u_3,v_3)-(u_3,v_2)-(u_3,v_1).$$

%$M_{(3,2)-(2,2)}$ together with the edges $(v_2,u_3), (v_1,u_3)$ and $(v_1,u_2)(v_1,u_3)$. 

Since the odd paths we have taken belong to different copies of $C$, they are pairwise edge-disjoint.

\end{proof}

\begin{corollary}
    $$\toi(K_t \Square K_s) = \begin{cases}
      t+s-1 &\text{ if }  t, s\geq 2 \text{ and } \max\{s,t\} \geq 4\\
      4 &\text{ if } s=t=3\\
      3 &\text{ if } s=3 \text{ and } t=2
    \end{cases}$$
\end{corollary}

\begin{proof}
    First, the degree of every vertex in $K_t \Square K_s$ is $t+s-2$, and thus, $\toi (K_t \Square K_s) \leq t+s-1$. When $s\geq 4 \text{ and } t\geq 2$, from Theorem \ref{thm:tslargeCart} we have the lower bound. The other cases come from the Theorems \ref{thm:t3s3Cart} and \ref{thm:t3s2Cart}.
\end{proof}

From the above three theorems we deduce the following, which amounts  the case of Cartesian product of Theorem~\ref{qus:toi}.  

\begin{corollary}
    Let $G$ and $H$ be two graphs with $\toi(G)=t\ge 2$ and $\toi(H)=s\ge 2$, respectively. Then, $\toi (G \Square H) \geq \toi(K_t \Square K_s)$.
\end{corollary}

\begin{proof}
    The maximum degree of a vertex in $K_t \Square K_s$ is $t+s-2$, and thus, $\toi (K_t \Square K_s) \leq t+s-1$. Therefore,  Theorem~\ref{thm:tslargeCart} gives $\toi (G \Square H) \geq \toi(K_t \Square K_s)$, when $t > 3$ or $s > 3$. So we may assume that both $s \leq 3$ and $t \leq 3$. If both $s$ and $t$ are $3$, then by Theorem~\ref{thm:t3s3Cart}, $\toi (G \Square H) \geq \toi(K_t \Square K_s)$. So we may assume that either $s=2$ or $t=2$. Without loss of generality, we can assume that $t=2.$ If $s$ is also $2$, then $K_t \Square K_s$ is a bipartite graph and we have $\toi (G \Square H) \geq \toi(K_t \Square K_s)=2$, since bipartite graphs exclude a $K_3$ as a totally odd strong immersion. The only case remaining is $s=3$ and $t=2$ which follows by Theorem~\ref{thm:t3s2Cart}.
\end{proof}

\begin{corollary}
    The cartesian product of two graphs is not a minimal counterexample to Conjeture~\ref{conj:church}.
\end{corollary}
\begin{proof}
    Let us assume that there are graphs $G$ and $H$ so that $G\Square H$ is a minimal counterexample to Conjecture \ref{conj:church}. Then, $G$ and $H$ satisfy that $\toi(G)\geq  \chi(G)$ and $\toi(H)\geq  \chi(H)$. Suppose that $\toi(G) \geq 3$ and $\toi(H)\geq 3$, since $\chi(G\Square H)\le \max\{\chi(G), \chi(H)\}\leq \chi(G)+\chi(H)-1 \leq \toi(G) + \toi(H)-1$, then from Theorem \ref{thm:tslargeCart}, we have $\chi(G\Square H)\le\toi(G) + \toi(H)-1 \le \toi(G\Square H)$, we obtain a contradiction with the fact that $G\Square H$ is a minimal counterexample to Conjecture \ref{conj:church}.

    Suppose that $\chi(G) \leq \toi(G)= 3$ and $\chi(H) \leq\toi(H)= 3$, then $\chi(G \Square H)\leq 3$, but from the Theorem \ref{thm:t3s3Cart} $\toi(G\Square H) \geq 4 \geq \chi(G \Square H)$, contradiction. Now, without loss of generality, suppose that $\chi(G) \leq \toi(G) = 3$ and $\chi(H) \leq\toi(H) = 2$, then clearly $\chi(G \Square H)\leq 3$ and from the Theorem \ref{thm:t3s2Cart} we get a contradiction. Finally, the case where $\chi(G) \leq \toi(G) =2$ and $\chi(H) \leq\toi(H) = 2$ is trivial as mentioned at the beginning of Section \ref{sec:cartesian}.
\end{proof}

\section{Acknowledgments}

 H. Echeverría is supported by ANID BECAS/DOCTORADO NACIONAL 21231147;  
  A. Jim\'enez is  supported by  ANID/Fondecyt Regular 1220071 and ANID-MILENIO-NCN2024-103;
  %S. Mishra is supported by ANID/Fondecyt Postdoctoral grant $3220618$ of Agencia Nacional de Investigati\'{o}n y Desarrollo (ANID), Chile;
  D.A. Quiroz is supported by ANID/Fondecyt Regular 1252197 and MATH-AMSUD MATH230035;
  M. Yépez is supported by ANID BECAS/DOCTORADO NACIONAL 21231444.
%

%\begin{comment}

%  \noindent E-mails:   henry.echeverria@postgrado.uv.cl (H. Echeverría),   andrea.jimenez@uv.cl (A. Jim\'enez),  suchismitamishra6@gmail.com (S. Mishra),  daniel.quiroz@uv.cl (D. Quiroz),  mauricio.yepez@postgrado.uv.cl (M. Yépez)

%  \end{comment}

\end{document}